\magnification 1200

\nopagenumbers
\headline={\ifnum\pageno=1 \hfill \else\hss{\tenrm--\folio--}\hss \fi}

\newcount\sectionnumber
\newcount\equationnumber
\newcount\thnumber
\newcount\refnumber

\def\ifundefined#1{\expandafter\ifx\csname#1\endcsname\relax}
\def\assignnumber#1#2{%
	\ifundefined{#1}\relax\else\message{#1 already defined}\fi
	\expandafter\xdef\csname#1\endcsname
 {\the\sectionnumber.\the#2}}%
\def\secassignnumber#1#2{%
  \ifundefined{#1}\relax\else\message{#1 already defined}\fi
  \expandafter\xdef\csname#1\endcsname{\the#2}}%
%
% macros on section numbers
%
\def\secname#1{\relax
  \global\advance\sectionnumber by 1
  \secassignnumber{S#1}\sectionnumber
  \csname S#1\endcsname}
\def\Sec#1 #2 {\vskip0pt plus.1\vsize\penalty-250\vskip0pt plus-.1\vsize
  \bigbreak\bigskip
  \equationnumber0\thnumber0
  \noindent{\bf \secname{#1}. #2}\par
  \nobreak\smallskip\noindent}
\def\sectag#1{\ifundefined{S#1}\message{S#1 undefined}{\sl #1}%
  \else\csname S#1\endcsname\fi}
\def\eq#1{\relax
  \global\advance\equationnumber by 1
  \assignnumber{EN#1}\equationnumber
  {\rm (\csname EN#1\endcsname)}}
\def\eqtag#1{\ifundefined{EN#1}\message{EN#1 undefined}{\sl (#1)}%
  \else(\csname EN#1\endcsname)\fi}
%
% macros on theorem numbers
%
\def\thname#1{\relax
  \global\advance\thnumber by 1
  \assignnumber{TH#1}\thnumber
  \csname TH#1\endcsname}
\def\thtag#1{\ifundefined{TH#1}\message{TH#1 undefined}{\sl #1}%
  \else\csname TH#1\endcsname\fi}
\def\Assumption#1 {\bLP{\bf Assumption \thname{#1}}\quad}
\def\Cor#1 {\bLP{\bf Corollary \thname{#1}}\quad}
\def\Def#1 {\bLP{\bf Definition \thname{#1}}\quad}
\def\Example#1 {\bLP{\bf Example \thname{#1}}\quad}
\def\Lemma#1 {\bLP{\bf Lemma \thname{#1}}\quad}
\def\Prop#1 {\bLP{\bf Proposition \thname{#1}}\quad}
\def\Remark#1 {\bLP{\bf Remark \thname{#1}}\quad}
\def\Theor#1 {\bLP{\bf Theorem \thname{#1}}\quad}
\def\Proof{\sLP{\bf Proof}\quad}
%
% macros on reference numbers
%
\def\refitem#1 #2\par{\ifundefined{REF#1}
\global\advance\refnumber by1%
\expandafter\xdef\csname REF#1\endcsname{\the\refnumber}%
\else\item{\ref{#1}}#2\sLP\fi}

\def\ref#1{\ifundefined{REF#1}\message{REF#1 is undefined}\else
  [\csname REF#1\endcsname]\fi}
\def\reff#1#2{\ifundefined{REF#1}\message{REF#1 is undefined}\else
  [\csname REF#1\endcsname, #2]\fi}
\def\Ref{\vskip0pt plus.1\vsize\penalty-250\vskip0pt plus-.1\vsize
  \bigbreak\bigskip\leftline{\bf References}\nobreak\smallskip
  \frenchspacing}
%
% other macros
%

\let\bPP=\bigbreak

\def\sLP{\smallbreak\noindent}

\def\bLP{\bigbreak\noindent}

\def\halmos{\hbox{\vrule height0.31cm width0.01cm\vbox{\hrule height
 0.01cm width0.3cm \vskip0.29cm \hrule height 0.01cm width0.3cm}\vrule
 height0.31cm width 0.01cm}}
\def\hhalmos{{\unskip\nobreak\hfil\penalty50
	\quad\vadjust{}\nobreak\hfil\halmos
	\parfillskip=0pt\finalhyphendemerits=0\par}}

\def\ga{\gamma}
\def\de{\delta}
\def\ep{\varepsilon}

\def\Ga{\Gamma}

\def\iy{\infty}
\def\const{{\rm const.}\,}
\def\Re{{\rm Re}\,}
\def\LHS{left-hand side}
\def\RHS{right-hand side}
\let\wt=\widetilde

\def\crs{\crcr\noalign{\smallskip}}
\def\crm{\crcr\noalign{\medskip}}
\overfullrule0pt
\font\titlefont=cmr10 at 17.28truept
\font\authorfont=cmr10 at 14truept
\font\addressfont=cmr10 at 10truept
\font\ttaddressfont=cmtt10 at 10truept

\centerline
{\titlefont
Identities of nonterminating series by Zeilberger's algorithm}
\bPP
\centerline{\authorfont Tom H. Koornwinder}
\vskip 2truecm
\noindent
{\sl Abstract:}\quad
This paper argues that automated proofs of identities for
non-terminating hypergeometric series are feasible by a combination
of Zeilberger's algorithm and asymptotic estimates.
For two analogues of Saalsch\"utz' summation
formula in the non-terminating case this is illustrated.

\bLP
{\sl Last modified:}\quad February 13, 1998
%ref.tex
\refitem{Erd53a}

\refitem{GasRah90}

\refitem{Ges95}

\refitem{Koe98}

\refitem{Koo90}

\refitem{Olv74}

\refitem{PeWZ}

\refitem{Zei97}

%intro.tex
\Sec{intro} {Introduction}
The Gosper algorithm and in particular the subsequent Zeilberger
algorithm and related WZ method have been extremely successful for an approach
by computer algebra to identities involving ($q$-)hypergeometric functions
(see the book by Petkov{\v s}ek, Wilf \& Zeilberger \ref{PeWZ}
and the references
given there).
Most of the current theory and all implementations of Zeilberger's
algorithm remain restricted to the case of identities for terminating
hypergeometric series, while many of the known identities in the literature
involve non-terminating hypergeometric identities. The book \ref{PeWZ}
discusses in Chapter 7 some nonterminating cases obtained from
terminating cases by the WZ method.
A very convincing demonstration that Zeilberger's algorithmic is
suitable for non-terminating cases is given by Gessel \reff{Ges95}{p.547}. He
demonstrates Gauss' summation formula for the Gaussian
hypergeometric series of argument 1 by means of a combination of
Zeilberger's algorithm and the asymptotic estimate
$$
{\Ga(a+k)\over\Ga(b+k)}\sim k^{a-b}\quad\hbox{as $k\to\iy$}
\eqno\eq{19}
$$
(see for instance Olver \reff{Olv74}{Ch.4, \S5.1} for a proof of \eqtag{19}).

I will demonstrate in this paper that the method extends to the
non-terminating generalization of Saalsch\"utz' summation formula for a
terminating Saalsch\"utzian ${}_3F_2$ hypergeometric series of argument 1.
This
nonterminating case has the form $A+B=C$, where $A$ is a non-terminating
Saalsch\"utzian ${}_3F_2(1)$, $B$ is a quotient of Gamma functions times
another
${}_3F_2(1)$ and $C$ is a quotient of Gamma functions. Without knowing the
formula explicitly, one can start with $A$, apply Zeilberger's algorithm
to the terms of $A$, and finally arrive by some limit transition
using \eqtag{19} at the desired formula $A+B=C$.
 
It is important to notice that all steps in the proof, both in the Gauss
case and in the non-terminating Saalsch\"utz case,
may be automated, including the application of the asymptotic formula
\eqtag{19}.
I strongly believe that this idea may lead to
an algorithmic approach to most identities for non-terminating
hypergeometric series in literature. For the $q$-analogues of such identities
an algorithmic approach should be feasible as well.

For an introduction to Zeilberger's and related algorithms I refer to
the book \ref{PeWZ}. For hypergeometric series see
Erd\'elyi \ref{Erd53a} and Gasper \& Rahman \ref{GasRah90}.
Although the emphasis in \ref{GasRah90} is on the $q$-case,
it also contains some information on the $q=1$ case.
Anyhow, many formulas for $q=1$ can be looked up from
the corresponding $q$-case in \ref{GasRah90},
by silently taking the formal limit for $q\uparrow1$, and by using
that
$$
\Ga_q(z):={(1-q)^{1-z}\,(q;q)_\iy\over(q^z;q)_\iy}\to\Ga(z)\quad
\hbox{as $q\uparrow 1$}
$$
(see \reff{Koo90}{Appendix B} and references given there).

The contents of this paper are as follows.
Section 2 discusses Gauss' summation formula following Gessel
\reff{Ges95}{p.547}.
Section 3 discusses the case of a non-terminating Saalsch\"utzian
${}_3F_2(1)$.
Section 4 gives three other identities related to the one derived in
section 3, and it is shown how these identitities follow analytically
from each other.
Finally, in section 5, we consider one of the other identities in
section 4 ($A+B=C$ with $C$ a Gamma quotient and $A$ and $B$ a
Gamma quotient times a non-terminating ${}_3F_2(1)$ with one of the
upper parameters equal to 1). As is demonstrated,
it can be derived by computer algebra
and asymptotics, similarly as for the case of section 3, but the
asymptotics is quite tricky (maybe interesting in its own right)
and not yet suitable to be automated.
%
%gauss.tex
%
\Sec{Gauss} {Gauss' summation formula}
Gauss' summation formula (Gauss, 1813; see for instance
\reff{Erd53a}{2.8(46) with proof in \S2.1.3}) is the non-terminating analogue
of the Chu-Vandermonde summation formula. It reads as follows:
$$
{}_2F_1\left[{a,b\atop c};1\right]:=
\sum_{k=0}^\iy {(a)_k\,(b)_k\over (c)_k\,k!}=
{\Ga(c)\,\Ga(c-a-b)\over \Ga(c-a)\,\Ga(c-b)}\,,
\eqno\eq{1}
$$
where we assume that $\Re(c-a-b)>0$ and $c\notin\{0,-1,-2,\ldots\}$
in order to ensure that the hypergeometric series on the \LHS\ is
well-defined and absolutely convergent.
Gessel \reff{Ges95}{\S7}
shows how to prove this identity by the WZ method,
as I will recapitulate now.

Formula \eqtag{1}, with $c$ replaced by $c+n$
($n\in\{0,1,2,\ldots\}$), can be written as
$$
\sum_{k=0}^\iy f(n,k)=s(n),
\eqno\eq{2}
$$
where
$$
f(n,k):={(a)_k\,(b)_k\over (c+n)_k\,k!}\,,\quad
s(n):={\Ga(c+n)\,\Ga(c+n-a-b)\over\Ga(c+n-a)\,\Ga(c+n-b)}\,.
\eqno\eq{6}
$$
The sum on the \LHS\ of \eqtag{2} still absolutely converges
(since $n\ge0$).
Formula \eqtag{2} can be rewritten as
$$
\sum_{k=0}^\iy F(n,k)=1,
\eqno\eq{3}
$$
where
$$
F(n,k):={f(n,k)\over s(n)}=
{\Ga(a+k)\,\Ga(b+k)\,\Ga(c+n-a)\,\Ga(c+n-b)\over
\Ga(1+k)\,\Ga(c+n+k)\,\Ga(a)\,\Ga(b)\,\Ga(c+n-a-b)}\,.
\eqno\eq{7}
$$
By \eqtag{19} we have
$$
F(n,k)\sim k^{a+b-c-1-n}\quad\hbox{as $k\to\iy$.}
$$
Thus the sum on the \LHS\ of \eqtag{3} converges absolutely
(as we already knew).
We want to prove \eqtag{3} by the WZ method.
Gosper's algorithm applied to $F(n+1,k)-F(n,k)$ or, equivalently, the WZ method
applied to $F(n,k)$ succeeds.
In Maple V4, for instance, call
$$
\eqalignno{
&\hbox{\tt read \`{}hsum.mpl\`{}: gosper(F(n+1,k)-F(n,k),k);}
\cr
\noalign{\hbox{or}}
&\hbox{\tt read \`{}hsum.mpl\`{}: WZcertificate(F(n,k),k,n);}
\cr}
$$
from Koepf's
package hsum.mpl \ref{Koe98}, or call
$$
\hbox{\tt read ekhad; ct(F(n,k),1,k,n,N);}
$$
from Zeilberger's package EKHAD \ref{Zei97}.
The resulting formula is:
$$
F(n+1,k)-F(n,k)=G(n,k+1)-G(n,k),
\eqno\eq{4}
$$
where
$$
G(n,k):={k\over c+n-a-b}\,F(n,k).
$$
Note that
$$
G(n,0)=0,\qquad
G(n,k)\sim {k^{a+b-c-n}\over c+n-a-b}\quad\hbox{as $k\to\iy$.}
\eqno\eq{5}
$$
It follows from \eqtag{4} and \eqtag{5} that
$$
\sum_{k=0}^K F(n+1,k)-\sum_{k=0}^K F(n,k)=
G(n,K+1)-G(n,0)=G(n,K+1)\to 0\quad\hbox{as $K\to\iy$.}
$$
Hence $\sum_{k=0}^\iy F(n,k)$ is independent of $n$. Thus
$$
\sum_{k=0}^\iy F(n,k)=
\lim_{n\to\iy}\sum_{k=0}^\iy F(n,k)=
\sum_{k=0}^\iy\Bigl(\lim_{n\to\iy} F(n,k)\Bigr)=\sum_{k=0}^\iy \de_{k,0}=1,
\eqno\eq{9}
$$
where the second equality is justified by dominated convergence.
Indeed, for some $C>0$ we have $1/|s(n)|\le C$ if $n\in\{0,1,2,\ldots\}$
(use \eqtag{6} and \eqtag{19}). Thus, by \eqtag{7} and \eqtag{6}
we obtain
$$
|F(n,k)|\le C\,|f(n,k)|\le\left|{(a)_k\,(b)_k\over (\Re c+n_0)_k\,k!}\right|
\eqno\eq{8}
$$
for $n\in\{n_0,n_0+1,\ldots\}$, where $n_0$ is such that $\Re c+n_0>0$.
Again by \eqtag{19}, the \RHS\ of \eqtag{8} summed over
$k=0,1,2,\ldots\;$ yields a convergent sum, thus justifying the
application of the dominated convergence theorem in \eqtag{9}.
%
%saalschutz.tex
%
\Sec{saalschutz} {Non-terminating analogue of Saalsch\"utz'
summation formula}
Saalsch\"utz' summation formula (which actually goes back to Pfaff, 1797)
for a terminating Saalsch\"utzian ${}_3F_2$ series with argument 1
reads as follows:
$$
\eqalign{
&{}_3F_2\left[{-m,b,c\atop e,-m+b+c-e+1};1\right]:=
\sum_{k=0}^m{(-m)_k\,(b)_k\,(c)_k\over (e)_k\,(-m+b+c-e+1)_k\,k!}
\crs
&\qquad={(e-b)_m(e-c)_m\over (e)_m(e-b-c)_m}\qquad
(m\in\{0,1,2,\ldots\},\quad e,e-b-c\notin\{0,-1,-2,\ldots\}).}
\eqno\eq{10}
$$
See for instance \reff{Erd53a}{4.4(3) with proof in \S2.1.5}.
Now replace $-m$ by $a$ on the \LHS\ of \eqtag{10}, where
$a$ is arbitrarily complex:
$$
{}_3F_2\left[{a,b,c\atop e,a+b+c-e+1};1\right]:=
\sum_{k=0}^m{(a)_k\,(b)_k\,(c)_k\over (e)_k\,(a+b+c-e+1)_k\,k!}\,,
\eqno\eq{11}
$$
and try to generalize identity \eqtag{10}
for this case. Such generalizations, with one additional
${}_3F_2(1)$ term, are known in literature (see for instance ...),
but we want to find a formula
of this type from scratch, just starting from \eqtag{11} and working
with the Zeilberger algorithm.

In \eqtag{11} replace $c$ by $c+n$, where $n\in\{0,1,2,\ldots\}$. Then
$$
{}_3F_2\left[{a,b,c+n\atop e,a+b+c-e+n+1};1\right]=
\sum_{k=0}^\iy f(n,k),
\eqno\eq{12}
$$
where
$$
f(n,k):={(a)_k\,(b)_k\,(c+n)_k\over (e)_k\,(a+b+c-e+1+n)_k\,k!}\,,
$$
so
$$
f(n,k)\sim
{\Ga(e)\,\Ga(a+b+c-e+1+n)\over\Ga(a)\,\Ga(b)\,\Ga(c+n)}\,k^{-2}\quad
\hbox{as $k\to\iy$}
$$
by \eqtag{19}. Thus we have absolute convergence in the sum in \eqtag{12}.
Zeilberger's algorithm applied to $f(n,k)$ succeeds.
In Maple V4, for instance, call
$$
\hbox{\tt read ekhad: ct(f(n,k),1,k,n,N);}
$$
from Zeilberger's package EKHAD \ref{Zei97}.
The resulting recurrence is:
$$
\eqalignno{
&(c-e+a+n+1)\,(b+c-e+n+1)\,f(n+1,k)
\cr
&\qquad-(c-e+n+1)\,(b+c-e+a+n+1)\,f(n,k)
=g(n,k+1)-g(n,k),&\eq{13}
\cr}
$$
where
$$
g(n,k):={(b+c-e+a+n+1)\,(e+k-1)\,k\over c+n}\,f(n,k).
$$
We can rewrite \eqtag{13} in the form
$$
F(n+1,k)-F(n,k)=G(n,k+1)-G(n,k),
\eqno\eq{14}
$$
where
$$
\eqalignno{
F(n,k):=&{\Ga(a+c-e+1+n)\,\Ga(b+c-e+1+n)\over\Ga(c-e+1+n)\,\Ga(a+b+c-e+1+n)}\,
f(n,k),
\crs
G(n,k):=&{k\,(e+k-1)\over (c+n)\,(c-e+n+1)}\,F(n,k).
\cr}
$$
Observe that
$$
F(n,k)\sim{\Ga(e)\,\Ga(a+c-e+1+n)\,\Ga(b+c-e+1+n)\over
\Ga(a)\,\Ga(b)\,\Ga(c+n)\,\Ga(c-e+1+n)}\,k^{-2}\quad\hbox{as $k\to\iy$},
$$
and
$$
G(n,0)=0,\quad
\lim_{k\to\iy}G(n,k)=
{\Ga(e)\,\Ga(a+c-e+1+n)\,\Ga(b+c-e+1+n)\over
\Ga(a)\,\Ga(b)\,\Ga(c+1+n)\,\Ga(c-e+2+n)}\,.
\eqno\eq{15}
$$
Put
$$
\eqalignno{
&S(n):=\sum_{k=0}^\iy F(n,k)
\crs
&={\Ga(a+c-e+1+n)\,\Ga(b+c-e+1+n)\over
\Ga(c-e+1+n)\,\Ga(a+b+c-e+1+n)}\,
{}_3F_2\left[{a,b,c+n\atop e,a+b+c-e+1+n};1\right].&\eq{17}
\cr}
$$
Then by \eqtag{14} and \eqtag{15} we obtain
$$
\eqalignno{
S(n+1)-S(n)&=\lim_{k\to\iy}G(n,k+1)-G(n,0)
\crs
&={\Ga(e)\,\Ga(a+c-e+1+n)\,\Ga(b+c-e+1+n)\over
\Ga(a)\,\Ga(b)\,\Ga(c+1+n)\,\Ga(c-e+2+n)}\,.&\eq{29}
\cr}
$$
Hence
$$
S(n)=S(0)+
{\Ga(e)\,\Ga(a+c-e+1)\,\Ga(b+c-e+1)\over\Ga(a)\,\Ga(b)\,\Ga(c+1)\,\Ga(c-e+2)}
\,\sum_{k=0}^{n-1}{(a+c-e+1)_k\,(b+c-e+1)_k\over(c+1)_k\,(c-e+2)_k}\,.
\eqno\eq{16}
$$
Note that
$$
{(a+c-e+1)_k\,(b+c-e+1)_k\over(c+1)_k\,(c-e+2)_k}\sim
{\Ga(c+1)\,\Ga(c-e+2)\over\Ga(a+c-e+1)\,\Ga(b+c-e+1)}\,k^{a+b-e-1}
\quad\hbox{as $k\to\iy$}.
$$
Therefore assume that $\Re(e-a-b)>0$. Then, for $n\to\iy$,
the \RHS\ of \eqtag{16} tends to
$$
\eqalignno{
&{\Ga(a+c-e+1)\,\Ga(b+c-e+1)\over\Ga(c-e+1)\,\Ga(a+b+c-e+1)}\,
{}_3F_2\left[{a,b,c\atop e,a+b+c-e+1};1\right]
\crm
&\qquad+{\Ga(e)\,\Ga(a+c-e+1)\,\Ga(b+c-e+1)\over
\Ga(a)\,\Ga(b)\,\Ga(c+1)\,\Ga(c-e+2)}\,
{}_3F_2\left[{a+c-e+1,b+c-e+1,1\atop c+1,c-e+2};1\right].
\cr}
$$
Here we substituted the \RHS\ of \eqtag{17} with $n=0$ for $S(0)$.
Next let $n\to\iy$ on the \LHS\ of \eqtag{16}.
From \eqtag{17} we obtain
$$
\lim_{n\to\iy} S(n)=
{}_2F_1\left[{a,b\atop e};1\right]=
{\Ga(e)\,\Ga(e-a-b)\over\Ga(e-a)\,\Ga(e-b)}\quad
(\Re(e-a-b)>0).
\eqno\eq{20}
$$
The second identity is Gauss' summation formula \eqtag{1}.
The first identity follows by applying \eqtag{19} to the limit of the
quotient of Gamma functions in \eqtag{17} combined with a formal
limit (to be justified in the Lemma below)
within the ${}_3F_2(1)$ in \eqtag{17}.
Thus the limit case of \eqtag{16} for $n\to\iy$ is the
identity we looked for:
$$
\eqalignno{
&{\Ga(a+c-e+1)\,\Ga(b+c-e+1)\over\Ga(c-e+1)\,\Ga(a+b+c-e+1)}\,
{}_3F_2\left[{a,b,c\atop e,a+b+c-e+1};1\right]
\crm
&\qquad+{\Ga(e)\,\Ga(a+c-e+1)\,\Ga(b+c-e+1)\over
\Ga(a)\,\Ga(b)\,\Ga(c+1)\,\Ga(c-e+2)}\,
{}_3F_2\left[{a+c-e+1,b+c-e+1,1\atop c+1,c-e+2};1\right]
\crm
&\qquad\qquad\qquad\qquad\qquad
={\Ga(e)\,\Ga(e-a-b)\over\Ga(e-a)\,\Ga(e-b)}\,,&\eq{18}
\cr}
$$
valid for $\Re(e-a-b)>0$.

Let us now give the promised Lemma which will justify the first
identity in \eqtag{20}.

\Lemma{21}
Let $\Re(e-a-b)>0$. Then
$$
\lim_{n\to\iy}{}_3F_2\left[{a,b,c+n\atop e,a+b+c-e+1+n};1\right]=
{}_2F_1\left[{a,b\atop e},1\right].
$$
\Proof
The limit formally holds because
$$
\lim_{n\to\iy}
{(a)_k\,(b)_k\,(c+n)_k\over (e)_k\,(a+b+c-e+1+n)_k\,k!}=
{(a)_k\,(b)_k\over (e)_k\,k!}\,.
$$
We will justify this limit by dominated convergence.
Take $\ep\in(0,1)$ such that also $\ep<\Re(e-a-b)$.
It follows from \eqtag{19} that
\goodbreak
$$
\eqalignno{
\left|{(a)_k\,(b)_k\,(c+n)_k\over (e)_k\,(a+b+c-e+1+n)_k\,k!}\right|
=\const\left|{\Ga(a+k)\,\Ga(b+k)\over\Ga(e+k)\,\Ga(1+k)}\right|\,
\qquad\qquad\qquad\qquad&
\crm
\times
\left|{\Ga(c+a+b-e+1+n)\over\Ga(c+n)}\right|\,
\left|{\Ga(c+n+k)\over\Ga(c+a+b-e+1+n+k)}\right|\qquad\qquad&
\crm
\le\const k^{\Re(a+b-e)-1}\,n^{\Re(a+b-e)+1}\,(n+k)^{\Re(e-a-b)-1}&
\crm
=\const \left({n\over n+k}\right)^{1-\ep}\,
k^{-1-\ep}\,(n^{-1}+k^{-1})^{\Re(e-a-b)-\ep}&
\crm
\le\cases{
\const k^{-\Re(e-a-b)-1}\le \const k^{-\ep-1}& if $k\le n$,
\crs
\const k^{-\ep-1}& if $k\ge n$.
\cr}&
\cr}
$$
Then the dominated convergence follows
because $\sum_{k=1}^\iy k^{-\ep-1}<\iy$.\hhalmos
%
%saalschutz2.tex
%
\Sec{saalschutz2} {Other non-terminating identities related to Saalsch\"utz'
summation formula}
Formula \eqtag{18}, which I derived in the previous section by a mixture
of computer algebra and asymptotic techniques, is possibly not in the
literature, but it can be derived from some formulas of similar nature
which are in the literature.

The best known non-terminating generalization of
Saalsch\"utz' formula \eqtag{10} is:
$$
\eqalignno{
&{}_3F_2\left[{a,b,c\atop e,a+b+c-e+1};1\right]
\crm
&\hskip 1truecm\relax
+{\Ga(e-1)\,\Ga(a-e+1)\,\Ga(b-e+1)\,\Ga(c-e+1)\,\Ga(a+b+c-e+1)
\over
\Ga(1-e)\,\Ga(a)\,\Ga(b)\,\Ga(c)\,\Ga(a+b+c-2e+2)}\,
\crm
&\hskip 2truecm\relax\times
{}_3F_2\left[{a-e+1,b-e+1,c-e+1\atop 2-e,a+b+c-2e+2};1\right]
\crm
&\hskip 3truecm\relax={\Ga(a-e+1)\,\Ga(b-e+1)\,\Ga(c-e+1)\,\Ga(a+b+c-e+1)\over
\Ga(1-e)\,\Ga(b+c-e+1)\,\Ga(a+c-e+1)\,\Ga(a+b-e+1)}\,,&\eq{22}
\cr}
$$
see \reff{GasRah90}{(II.24)}.
Since both ${}_3F_2(1)$ series are Saalsch\"utzian, there are no
conditions on the parameters necessary for convergence.

The second formula we will use involves two non-terminating
${}_3F_2(1)$ series with
one upper parameter equal to 1, like the second ${}_3F_2(1)$ in formula
\eqtag{18}. The formula reads as follows.
$$
\eqalignno{
{}_3F_2\left[{a,b,1\atop d,e};1\right]+
{1-e\over e-b-1}\,{}_3F_2\left[{d-a,b,1\atop d,b-e+2};1\right]
\hskip 4truecm\relax&
\crm
=
{\Ga(d)\,\Ga(e)\,\Ga(a-e+1)\,\Ga(b-e+1)\,\Ga(d+e-a-b-1)\over
\Ga(a)\,\Ga(b)\,\Ga(d-a)\,\Ga(d-b)}\qquad&&\eq{23}
\cr\noalign{\vskip 6pt}
(\Re(d+e-a-b-1)>0,\;\Re(a-e+1)>0).&
\cr}
$$
It can be derived by specialization of \reff{GasRah90}{(3.3.1)}
combined with Gauss' summation formula \eqtag{1}.

The third formula transforms a nonterminating Saalsch\"utzian
${}_3F_2(1)$ into a nonterminating ${}_3F_2(1)$ with upper parameter 1:
$$
\eqalignno{
{}_3F_2\left[{a,b,c\atop e,a+b+c-e+1};1\right]=
{\Ga(e)\,\Ga(a+b+c-e+1)\over\Ga(a)\,\Ga(b+1)\,\Ga(c+1)}\hskip 4truecm\relax&
\crm
\times{}_3F_2\left[{e-a,b+c-e+1,1\atop b+1,c+1};1\right]\qquad
(\Re a>0).&&\eq{24}
\cr}
$$
It follows by specialization of \reff{GasRah90}{(3.1.2)}.

Now replace in \eqtag{24} $a,b,c,e$ by
$b-e+1,a-e+1,c-e+1,2-e$. Next substitute in the \RHS\ of \eqtag{24} (with
parameters replaced as above) formula \eqtag{23} with
$a,b,d,e$ replaced by $1-b,a+c-e+1,a-e+2,c-e+2$.
This yields
$$
\eqalignno{
{}_3F_2\left[{b-e+1,a-e+1,c-e+1\atop 2-e,a+b+c-2e+2};1\right]=
-{\Ga(2-e)\,\Ga(a+b+c-2e+2)\over c\,\Ga(b-e+1)\,\Ga(a-e+1)\,\Ga(c-e+2)}&
\crm
\times{}_3F_2\left[{1+b+c-e,1+a+c-e,1\atop 2+c-e,c+1};1\right]
\hskip 3truecm\relax&
\crm
+
{\Ga(e-a-b)\,\Ga(c)\,\Ga(2-e)\,\Ga(a+b+c-2e+2)\over
\Ga(1-b)\,\Ga(1+b+c-e)\,\Ga(1-a)\,\Ga(a+c-e+1)}
\hskip 1truecm\relax&&\eq{25}
\cr\noalign{\vskip 6pt}
(\Re(1+b-e)>0,\;\Re(e-a-b)>0).&
\cr}
$$
Now formula \eqtag{18} follows from \eqtag{22} by substituting \eqtag{25}
for the second ${}_3F_2(1)$ in \eqtag{22}, and by using that
$$
\Ga(z)\,\Ga(1-z)={\pi\over\sin(\pi z)}\,.
$$

So we have seen that formulas \eqtag{22}, \eqtag{23} and \eqtag{24}
imply formula \eqtag{18} via formula \eqtag{25}.
Conversely, formulas \eqtag{18}, \eqtag{23} and \eqtag{24} imply formula
\eqtag{22} via formula \eqtag{25}. In the next section we will derive
formulas \eqtag{23} and \eqtag{24} by methods of computer algebra and
asymptotics.
%
%saalschutz3.tex
%
\Sec{saalschutz3} {Further proofs of identities of
non-terminating series by computer algebra and asymptotics}
Let us try to prove formula \eqtag{23} in a similar way as formula \eqtag{18}.
Put
$$
f(n,k):={(a)_k\,(b+n)_k\over(d+n)_k\,(e)_k}\,.
$$
By Zeilberger's algorithm we find that
$$
\eqalignno{
(a-n-d)\,(n&+b)\,f(n+1,k)-(n+b-e+1)\,(n+d)\,f(n,k)=g(n,k+1)-g(n,k),\quad
\crm
&\qquad\qquad\qquad\qquad\qquad\qquad\hbox{where}\quad
g(n,k):=-(n+d)\,(e+k-1)\,f(n,k).
\cr
\noalign{\hbox{Hence}}
&\qquad\qquad F(n+1,k)-F(n,k)=G(n,k+1)-G(n,k),
\cr
\noalign{\hbox{where}}
&F(n,k):=
{(b)_n\,(d-a)_n\over (d)_n\,(b-e+1)_n}\,
{(a)_k\,(b+n)_k\over (d+n)_k\,(e)_k}\,,
\crm
&G(n,k):=
{(e+k-1)\,(a)_k\,(b+n)_k\,(b)_n\,(d-a)_n\over
(e-n-b-1)\,(d+n)_k\,(e)_k\,(d)_n\,(b-e+1)_n}\,.
\cr}
$$
Assume
$$
\Re(d+e-a-b-1)>0,\quad
\Re(a-e+1)>0.
\eqno\eq{36}
$$
Then
$$
S(n):=\sum_{k=0}^\iy F(n,k)=
{(b)_n\,(d-a)_n\over (d)_n\,(b-e+1)_n}\,
{}_3F_2\left[{a,b+n,1\atop d+n,e};1\right]
\eqno\eq{28}
$$
has absolutely convergent series and
$$
\eqalignno{
S(n+1)-S(n)&=\lim_{k\to\iy}G(n,k+1)-G(n,0)=-G(n,0)
\cr
&=-{(e-1)\,(d-a)_n\,(b)_n\over(e-b-1)\,(d)_n\,(b-e+2)_n}\,.&\eq{30}
\cr}
$$
Hence
$$
\eqalignno{
S(0)=S(\iy)+{e-1\over e-b-1}\,\sum_{n=0}^\iy
{(d-a)_n\,(b)_n\over(d)_n\,(b-e+2)_n}\,,\qquad&
\crs
{\rm where}\quad S(\iy):=\lim_{n\to\iy} S(n).&
\cr}
$$
Thus
$$
{}_3F_2\left[{a,b,1\atop d,e};1\right]=
{e-1\over e-b-1}\,{}_3F_2\left[{d-a,b,1\atop d,b-e+2};1\right]+S(\iy).
$$
So formula \eqtag{23} will be proved if we can show that
$$
\eqalignno{
&\lim_{n\to\iy}{(b)_n\,(d-a)_n\over (d)_n\,(b-e+1)_n}\,
{}_3F_2\left[{a,b+n,1\atop d+n,e};1\right]
\crm
&\hskip 2truecm\relax
={\Ga(d)\,\Ga(e)\,\Ga(a-e+1)\,\Ga(b-e+1)\,\Ga(d+e-a-b-1)\over
\Ga(a)\,\Ga(b)\,\Ga(d-a)\,\Ga(d-b)}\,.&\eq{26}
\cr}
$$
This limit result is not evident. A formal limit would yield
$$
\left(\lim_{n\to\iy} n^{e-a-1}\right)
{}_2F_1\left[{a,1\atop e};1\right]=0\cdot\iy\quad
\hbox{since $\Re(a-e+1)>0$.}
$$
We will give a proof of \eqtag{26} at the end of this section.

Let us first turn to
formula \eqtag{24}, and rewrite it in the following equivalent form.
$$
\eqalignno{
&{(b)_n\,(d-a)_n\over (d)_n\,(b-e+1)_n}\,
{}_3F_2\left[{a,b+n,1\atop d+n,e};1\right]
\crm
&\qquad=
{\Gamma(d)\,\Gamma(e)\,\Gamma(b-e+1)\,\Gamma(d+e-a-b-1)\over
\Gamma(b)\,\Gamma(d-a)\,\Gamma(d+e-b-1)}\,
{\Gamma(b+n)\,\Gamma(d-a+n)\over
\Gamma(b-e+1+n)\,\Gamma(d+e-a-1+n)}\,
\crm
&\qquad\quad\times
{}_3F_2\left[{d+n-1,d+e-a-b-1,e-1\atop d+e-a+n-1,d+e-b-1};1\right]
\quad(\Re(d+e-a-b-1)>0).&\eq{27}
\cr}
$$
Observe that we can take limits on the \RHS\ of \eqtag{27} for $n\to\iy$
(see \eqtag{20}), and thus arrive at \eqtag{26}.
However, I only want to use \eqtag{27} if I can derive it by methods
of computer algebra and asymptotics. Let us try.

Write \eqtag{27} more compactly as
$$
S(n)=C\,\wt S(n),
\eqno\eq{31}
$$
where $S(n)$ is given by \eqtag{28},
$$
C:={\Gamma(d)\,\Gamma(e)\,\Gamma(b-e+1)\,\Gamma(d+e-a-b-1)\over
\Gamma(b)\,\Gamma(d-a)\,\Gamma(d+e-b-1)}\,,
$$
and
$$
\wt S(n):=
{\Gamma(b+n)\,\Gamma(d-a+n)\over
\Gamma(b-e+1+n)\,\Gamma(d+e-a-1+n)}\,
{}_3F_2\left[{d+n-1,d+e-a-b-1,e-1\atop d+e-a+n-1,d+e-b-1};1\right].
$$
Note that $\wt S(n)$ is just \eqtag{17} with
$a,b,c,e$ replaced by $e-1,d+e-a-b-1,d-1,d+e-b-1$.
Thus \eqtag{29} yields that
$$
\wt S(n+1)-\wt S(n)=
{1\over C}\,{(e-1)\,(d-a)_n\,(b)_n\over(b-e+1)\,(d)_n\,(b-e+2)_n}\,.
$$
Comparison with \eqtag{30} yields that
$S(n)-C\,\wt S(n)$ is independent of $n$. Hence
$$
S(n)-C\,\wt S(n)=S(\iy)-C\,\wt S(\iy).
$$
Thus \eqtag{31} holds iff \eqtag{26} holds.
So we have not made real progress. Formula \eqtag{26} will follow from
\eqtag{31}, but for the proof of \eqtag{31} we need \eqtag{26}.
Let us therefore give an independent proof of \eqtag{26}.

\bLP
{\bf Proof of \eqtag{26}}\quad
Observe that
$$
\eqalignno{
&{(b)_n\,(d-a)_n\over (d)_n\,(b-e+1)_n}\,
{}_3F_2\left[{a,b+n,1\atop d+n,e};1\right]
\crm
&\qquad={\Ga(d)\,\Ga(e)\,\Ga(b-e+1)\over\Ga(a)\,\Ga(b)\,\Ga(d-a)}\,
\sum_{k=0}^\iy
{\Ga(a+k)\over\Ga(e+k)}\,{\Ga(b+n+k)\over\Ga(d+n+k)}\,
{\Ga(d-a+n)\over\Ga(b-e+1+n)}\,.&\eq{32}
\cr}
$$
We will first work formally, and consider the sum obtained by
replacing in the above infinite sum
the gamma quotients by their asymptotic estimates
(for big $n$ and $k$) as given by \eqtag{19}. This yields
$$
\eqalignno{
&\sum_{k=0}^\iy
{\Ga(a+k)\over\Ga(e+k)}\,{\Ga(b+n+k)\over\Ga(d+n+k)}\,
{\Ga(d-a+n)\over\Ga(b-e+1+n)}&\eq{37}
\cr
&\sim\sum_{k=0}^\iy
(k+1)^{a-e}\,(n+k+2)^{b-d}\,(n+1)^{d+e-a-b-1}&\eq{33}
\cr
&=
{1\over n+1}\sum_{k=0}^\iy
\left(1+{k+1\over n+1}\right)^{b-d}\,
\left({k+1\over n+1}\right)^{a-e}&\eq{34}
\cr
&\sim
\int_{t=0}^\iy (1+t)^{b-d}\,t^{a-e}\,dt\quad\hbox{as $n\to\iy$}&\eq{35}
\cr\crs
&={\Ga(a-e+1)\,\Ga(d+e-a-b-1)\over\Ga(d-b)}\,.
\cr}
$$
Here we used that the integral in \eqtag{35}, which converges because of
assumptions \eqtag{36}, is approximated by its Riemann sum \eqtag{34}
as $n\to\iy$.
Thus formally the \RHS\ of \eqtag{32} indeed equals the \RHS\ of
\eqtag{26}.

We will justify this formal proof by Lebesgue's dominated convergence
theorem. Rewrite \eqtag{37} as
$$
\int_0^\iy \ga_n(t)\,h_n(t)\,dt,
\eqno\eq{38}
$$
where, for ${k\over n+1}<t\le {k+1\over n+1}$,
$$
h_n(t):=\left(1+{k+1\over n+1}\right)^{b-d}\,
\left({k+1\over n+1}\right)^{a-e}
$$
and
$$
\ga_n(t):=
{\Ga(a+k)\over\Ga(e+k)\,(k+1)^{a-e}}\,
{\Ga(b+n+k)\over\Ga(d+n+k)\,(n+k+2)^{b-d}}\,
{\Ga(d-a+n)\over\Ga(b-e+1+n)\,(n+1)^{d+e-a-b-1}}\,.
$$
Then
$$
|\ga_n(t)|\le\const\quad{\rm and}\quad
|h_n(t)|\le \const (1+t)^{\Re(b-d)}\,t^{\Re(a-e)}.
$$
Furthermore,
$$
\lim_{n\to\iy}\ga_n(t)=1\quad{\rm and}\quad
\lim_{n\to\iy}h_n(t)=(1+t)^{b-d}\,t^{a-e}.
$$
It follows by dominated convergence that the integral \eqtag{38} tends
to the integral \eqtag{35} as $n\to\iy$.\hhalmos
\Ref
%ref.tex
%
\refitem{Erd53a}
A. Erd\'elyi e.a.,
{\sl Higher transcendental functions, Vol. I},
McGraw-Hill, 1953. Reprinted in 1981 by R.E. Krieger.

\refitem{GasRah90}
G. Gasper \& M. Rahman,
{\sl Basic hypergeometric series},
Cambridge University Press, 1990.

\refitem{Ges95}
I. M. Gessel,
{\sl Finding identities with the WZ method},
J. Symb. Comp. 20 (1995), 537--566.

\refitem{Koe98}
W. Koepf,
{\sl Maple package hsum.mpl},
Version 1.0, February 1, 1998; obtainable by email from
{\tt koepf@imn.htwk-leipzig.de}.

\refitem{Koo90}
T. H. Koornwinder,
{\sl Jacobi functions as limit cases of $q$-ultraspherical polynomials},
J. Math. Anal. Appl. 148 (1990), 44--54.

\refitem{Olv74}
F. W. J. Olver,
{\sl Asymptotics and special functions},
Academic Press, New York, 1974.
Reprinted in 1997 by A. K. Peters.

\refitem{PeWZ}
M. Petkov{\v s}ek, H. S. Wilf and D. Zeilberger, 
$A=B$,\quad
A. K. Peters, 1996.

\refitem{Zei97}
D. Zeilberger,
{\sl Maple package EKHAD},
Version of March 27, 1997; obtainable from URL
{\tt http://www.math.temple.edu/\char'176 zeilberg}.

\nonfrenchspacing
\vskip 2 truecm
{\addressfont
\obeylines\parindent 5truecm
University of Amsterdam, Korteweg-de Vries Institute for Mathematics
Plantage Muidergracht 24, 1018 TV Amsterdam, The Netherlands
email {\ttaddressfont thk@wins.uva.nl}}
\bye